\documentclass[french,11pt,a4paper]{article}
\usepackage[utf8]{inputenc}

\usepackage[francais]{babel}
\usepackage[T1]{fontenc}

\usepackage{amsmath,amsthm}
\newtheorem*{theorem}{Théorème}

\newtheorem*{cor}{Corollaire}
\newtheorem*{rmq}{Remarque}
\usepackage{amsfonts}
\usepackage{amssymb}
\usepackage{babelbib}

\newcommand{\GL}{\mathrm{GL}}
\newcommand{\SL}{\mathrm{SL}}
\title{Sur les paquets d'Arthur aux places réelles, translation}
\date{}

  \author{Colette Moeglin\\
 CNRS, Institut de Mathématiques de Jussieu \\
 colette.moeglin@imj-prg.fr
\and
David Renard  \\Centre de Mathématiques
  Laurent Schwartz,  Ecole Polytechnique\\
david.renard@polytechnique.edu}
\date{}

\begin{document}
\maketitle

\section{Introduction}Pour décrire les formes automorphes de carré intégrable d'un groupe réductif $G$ ou plutôt les représentations des points adéliques de ce  groupe qui se réalisent dans cet espace, Arthur a introduit des paramètres locaux. Ces paramètres locaux  sont des morphismes de $W'_F\times \mathrm{SL}_2(\mathbb{C})$ dans le $L$-groupe de $G$ où $F$ est le corps local et $W'_F$ son groupe de Weil-Deligne. Ces paramètres en plus des propriétés décrites ci-dessous en  \ref{parametre} doivent avoir une compatibilité globale; cette compatibilité ne compte pas pour nous ici et on peut dire que pour le moment elle n'est décrite qu'artificiellement pour les groupes classiques et leurs variantes. Il faudra attendre la description des groupes tannakiens associés à des corps de nombres pour l'exprimer plus généralement. 

Soit $\psi$ un tel paramètre, on note $A(\psi)$ le groupe des composantes du centralisateur de $\psi$ dans la composante neutre de $^LG$, en fait plus générale\-ment, il faut passer à un rev\^etement (cf. \cite{book} chapitre 9) mais qui n'intervient pas pour les groupes quasi-déployés ni pour les groupes classiques m\^eme non quasi-déployés.

Arthur suggère qu'à tout tel paramètre, $\psi$, sur le corps local $F$, est associé une représentation de $G(F)\times A(\psi)$, unitaire, semi-simple par définition et de longueur finie. Pour les groupes classiques que l'on va considérer Arthur a  résolu tous ces problèmes locaux (cf. \ref{defpaquet} où l'on rappelle les résultats de \cite{book}). On  note $\pi^A(\psi)$  la représentation  associée à $\psi$.

Cette représentation  doit être compatible à l'endoscopie, à l'endoscopie tordue et à l'induction parabolique. On rappelle précisément en \ref{endoscopie} et \ref{endtord} ce que cela veut dire. Il faut remarquer que la définition de $\pi^A(\psi)$ dépend du choix des facteurs de transfert géométriques. Pour les groupes classiques, il n'y a aucun problème de choix pour ces facteurs de transfert, on suit Kottwitz et Shelstad. 

Toutefois, Arthur n'a défini $\pi^A(\psi)$ que pour les groupes classiques quasi-déployés. Pour un groupe classique non nécessairement quasi-déployé, les compatibilités requises à l'endoscopie et  l'endoscopie tordue déterminent uniquement $\pi^A(\psi)$.
C'est le fait que $\pi^A(\psi)$ soit bien une combinaison linéaire à coefficients positifs de représentations irréductibles de $G(F)\times A(\psi)$ (c'est à dire une représentation semi-simple de ce produit de groupes) qui n'est pas clair. C'est la conjecture 9.4.2 de \cite{book}. 

Via la stabilisation des formules de traces disponibles, donc des méthodes globales l'existence des représentations $\pi^A(\psi)$ est démontrée pour les groupes classiques quasi-déployés en \cite{book} 1.5.1 et 2.2.1 (b). On généralisera aux cas des groupes classiques non quasi-déployés dans un article ultérieur en copiant les méthodes de \cite{book}. Et en donnant une description explicite des représentations $\pi^A(\psi)$, on finira par démontrer que $\pi^A(\psi)$ est bien une représentation de $G(F)\times A(\psi)$.

Remarquons que pour les corps archimédiens \cite{ABV} donne aussi une définition de $\pi^A(\psi)$ beaucoup plus géométrique, qui ne suppose pas que $G$ est un groupe classique. La compatibilité des deux définitions pour les groupes classiques est en bonne voie avec les travaux de Mezo et ses collaborateurs, et on l'a montrée  dans \cite{MR} pour le corps des complexes.

Cet article poursuit sans l'atteindre le but suivant: démontrer que, pour les groupes classiques,  la représentation $\pi^A(\psi)$ est sans multiplicité, c'est-à-dire que c'est une somme de représentations irréductibles toutes distinctes en tant que représentations de $G(F)\times A(\psi)$ et même plus à savoir que si l'on décompose $\pi^A(\psi)$ dans la base des représentations  irréductibles de $G(F)$ alors les coefficients sont des représentations de dimension un (c'est-à-dire des caractères) de $A(\psi)$. Ce résultat est connu dans le cas quasi-déployé pour tout corps local sauf $F=\mathbb{R}$; on enlèvera sans difficulté l'hypothèse quasi-déployé pour les corps p-adiques car tout se ramène au cas des séries discrètes qui est assez facile (cf. \cite{note} section 3).

 Ici on suppose donc que le corps local est  $\mathbb{R}$. Ainsi  $W'_F$ n'est autre que le groupe de Weil $W_\mathbb{R}$ de $\mathbb{R}$. Expliquons un peu plus notre point de vue et ce que l'on obtient comme résultats.

On voit $\psi$ comme une représentation de $W_\mathbb{R}\times \mathrm{SL}(2,\mathbb{C})$ dans $\GL(\mathrm{dim}\,\psi,\mathbb{C})$ où $dim\,\psi$ est la dimension de la représen\-tation naturelle de $ ^LG$; cette repré\-sen\-tation est autoduale. On dit que $\psi$ est de bonne parité si toutes les sous-représentations irréductibles de $\psi$ sont autoduales et symplectiques si $^LG$ est un groupe symplectique, orthogonales sinon. 

En principe le cas général se ramène au cas de bonne parité par une induction irréductible; c'est un point technique que l'on n'aborde pas dans cet article. 

On s'intéresse donc au cas de bonne parité. Soit $\rho$ une représentation irréductible de  $W_\mathbb{R}$ que l'on suppose autoduale; la restriction de $\rho$ au sous-groupe $\mathbb{C}^*$ de $W_\mathbb{R}$ est soit le caractère nécessairement trivial soit la somme de deux caractères unitaires. Dans le premier cas, on pose $t_\rho=0$ et dans le deuxième cas, on note $t_\rho$ le demi-entier positif qui définit l'un des deux caractères. Soit $\psi$ un morphisme de bonne parité, on a une décompostion de $\psi$ en somme de représentations irréductibles (pour tout entier positif $a$ on note $R[a]$ la représentation irréductible de $\mathrm{SL}(2,\mathbb{C})$ de dimension $a$) $$
\psi=\bigoplus_{i\in [1,\ell]}\rho_i\otimes R[a_i],\eqno(*)
$$où les $a_i$ sont des entiers et les $\rho_i$ sont des représentations irréductibles autoduales de $W_\mathbb{R}$. On suppose ce qui est loisible que $t_{\rho_1}\geq \cdots\geq t_{\rho_\ell}$. On dit que $\psi$ est très régulier si dans la décomposition ci-dessus, il existe un entier $i_0\in [0,\ell]$ tel que pour tout $i>i_0$, $t_{\rho_i}=0$ et $t_{\rho_1}>> \cdots >> t_{\rho_{i_0}}>>0$. Si $i_0=0$, on dit que l'on est dans le cas spécial unipotent (en suivant Barbasch Vogan et \cite{pourhowe}).

Quand  cette condition de régularité  est satisfaite les méthodes de \cite{AMR} nous permettent dans \cite{recent} de calculer $\pi^A(\psi)$. On y montre  que cette représentation s'obtient par induction cohomologique explicite de représen\-ta\-tions faiblement unipoten\-tes au sens de Barbasch-Vogan. Un tel résultat est une généralisation de \cite{AJ}. Le cas spécial unipotent est compris pour les groupes quasi-déployés dans \cite{pourhowe} et est généralisé sans l'hypothèse quasi-déployé (cf. \cite{note} section 4). Le but de cet article est donc de ramener le cas général au cas très régulier.

Pour cela on va utiliser la translation. Fixons $G$ et $\psi$, on a donc la dimension de la représentation définie par $\psi$. Soit $\psi_+$  une autre représentation de $W_\mathbb{R}\times \SL(2,\mathbb{C})$, très régulière de même dimension que  $\psi$. On dit que $\psi_+$ domine $\psi$ si $\psi_+$ a la m\^eme décomposition (*) que $\psi$ mais où  les $\rho_i$ tels que $t_{\rho_i}\neq 0$ sont changés  en des $\rho'_i$ correspondant tels que $t_{\rho'_1}>> \cdots >>t_{\rho'_{i_0}}>>0$ où  $i_0= sup \{i\in [1,\ell]; t_{\rho_i}\neq 0$.

Dans cet article on montre que $\pi^A(\psi)$ s'obtient par translation à partir de $\pi^A(\psi_+)$ où $\psi_+$ est très régulier et domine $\psi$.  Si $i_0=0$, on est dans le cas des paramètres spéciaux unipotents traités en \cite{pourhowe} et \cite{note} et il n'y a pas de translation à faire. 

Il est facile de montrer ce résultat en prenant la translation déjà utilisée dans \cite{HS}, c'est-à-dire que l'on fixe $\psi,\psi_+$,  pour tout $i\in [1,i_0]$ on pose 
$T_i:= t_{\rho'_i}-t_{\rho_i}$. Ces nombres $T_i$ sont nécessairement des entiers positifs grâce à l'hypothèse de bonne parité.  Et on note $\lambda$ le poids,
$$(T_1, \cdots, T_1, \cdots , T_{i_0}, \cdots, T_{i_0}, 0, \cdots, 0),
$$où chaque $T_i$ vient $a_i$ fois. On considère le caractère d'un tore maximal de $G(\mathbb{C})$ égal à 
$
\sum_w e^{-w\lambda}$ où $w$ parcourt le groupe de Weyl de $G(\mathbb{C})$. On suit en cela \cite{HS} qui remarque que ce caractère du tore est la restriction au tore du caractère d'une combinaison linéaire de représentations de dimension finie de $G(\mathbb{C})$ et la translation se fait en tensorisant par cette combinaison linéaire de représentation de dimension finie restreinte à $G(\mathbb{R})$ et en projetant le résultat sur le caractère infinitésimal défini par $\psi$. Cette opération se fait dans le groupe de Grothendieck des modules de Harish-Chandra de $G(\mathbb{R})$.

On améliore ce résultat en utilisant la description de $\pi^A(\psi_+)$ que nous venons d'obtenir en \cite{recent}. On explique plus précisément dans le texte ci-dessous. Sous cette hypothèse, on montre que l'on peut aussi  utiliser le \og vrai\fg{} fonc\-teur de translation obtenu grâce à la représentation de $G(\mathbb{C})$ de plus bas poids $-\lambda$. C'est beaucoup mieux car ce foncteur de translation est lui défini directement dans la catégorie des modules de Harish-Chandra de longueur finie et est un foncteur exact.

 Et le résultat prouvé ici donnera donc une description de $\pi^A(\psi)$ en général.  Mais évidemment le foncteur de translation n'est pas simple à comprendre et les multiplicités un cherchées ne se déduisent pas aisément (c'est le moins que l'on puisse dire) du cas où $\psi$ est très régulier. 
 
 Pour expliquer un peu plus cette propriété de multiplicité un, soit $\epsilon$ un caractère de $A(\psi)$ et notons $\pi^A(\psi,\epsilon)$ l'unique représentation virtuelle ou 0 de $G(\mathbb{R})$ tel que (cf. \ref{defpaquet} ci-dessous pour plus de précision)
$$
\pi^A(\psi)=\sum_\epsilon \pi^A(\psi,\epsilon)\otimes \epsilon.
$$
Alors on montre que $\pi^A(\psi,\epsilon)$ se déduit de $\pi^A(\psi', \epsilon)$ pour $\psi'$ très régulier et convenablement relié à $\psi$ par un foncteur de translation, le foncteur de translation usuel de Zuckerman. Cela a en particulier comme conséquence que $\pi^A(\psi,\epsilon)$ est une combinaison linéaire à coefficients positifs ou nuls de représentations irréductibles unitaires, puisque l'on connaît ce résultat dans le cas très régulier. 
La propriété de multiplicité un cherchée pour $\psi$,   dit que les représentations $\pi^A(\psi,\epsilon)$ sont sans multiplicité et disjointes quand $\epsilon$ varie. Aucune de ces deux propriétés n'est clairement respectée par translation car la translation que l'on doit utiliser traverse des murs. Et bien que cette propriété de multiplicité un est clair dans le cas très régulier, on ne l'obtient pas pour $\psi$ général. On espère toutefois que la description obtenue ici aidera pour comprendre $\pi^A(\psi)$ au moins qualitativement.

\ 

La méthode est totalement élémentaire, consistant à utiliser le fait, qu'à condition de définir correctement les choses, la translation commute à l'endoscopie et à l'endoscopie tordue. C'est ce qui est expliqué en \ref{endoscopie}, \ref{endtord} et \ref{translation}. On montre  que ces propriétés entraînent aussi que si le paramètre $\psi$ est unipotent au sens de Barbasch-Vogan, c'est-à-dire que la restriction de $\psi$ au sous-groupe $\mathbb{C}^*$ de $W_\mathbb{R}$ est triviale, alors $\pi^A(\psi)$ restreint en une représentation de $G(\mathbb{R})$ est faiblement unipotente (\cite{KV} 12.3) ce qui était attendu par Barbasch-Vogan. Et cela nous sert pour conclure la démonstration du résultat que l'on vient d'expliquer, c'est-à-dire passer de la translation de \cite{HS} à la translation usuelle.

\thanks{Colette Moeglin remercie la fondation Simons et les organisateurs du congrès ayant eu lieu au château d'Elmau. David Renard a bénéficié d'une aide de  l'agence nationale de la recherche 
ANR-13-BS01-0012 FERPLAY. Les deux auteurs remercient le référé pour sa relecture soigneuse et constructive du texte.}

\section{Paquets d'Arthur}

\subsection{Paramètre d'Arthur}\label{parametre}
Un paramètre d'Arthur pour un groupe algébrique réductif $ G $ défini sur $ \mathbb{R} $ est un morphisme de groupes
\[\psi: W_\mathbb{R}\times \mathrm{SL}_2(\mathbb{C})\rightarrow \, ^LG \]
tel que 

(i) la restriction de $ \psi $ à $W_\mathbb{R}$ est un paramètre de Langlands tempéré,

(ii) la restriction de $\psi$ à $\mathrm{SL}_2(\mathbb{C})$ est un morphisme algébrique.

\noindent
A tout paramètre d'Arthur, on associe  un caractère infinitésimal: on considère la restriction de $\psi$ à $\mathbb{C}^*$ où l'on envoie $\mathbb{C}^*$ dans $W_\mathbb{R}\times \mathrm{SL}_2(\mathbb{C})$ par

$$ z\in \mathbb{C}^*\mapsto \biggl (z,\left(\begin{matrix}
  (z\overline{z})^{1/2}&0 \\ 0 &(z\overline{z})^{-1/2}\end{matrix}\right)\biggr).
$$
On construit ainsi un morphisme de $\mathbb{C}^*$ dans le groupe dual d'un tore maximal $T$ de $G$, donc une représentation de dimension un de $T(\mathbb{C})$. En identifiant cette représentation à sa différentielle, ce que l'on peut faire puisque l'on est sur $\mathbb{C}$, on obtient une classe de conjugaison d'éléments semi-simples de ${\mathrm{Lie}}(G)\otimes_\mathbb{R}\mathbb{C}$ et c'est ce que l'on appelle un caractère infinitésimal. Plus loin on  verra plutôt les caractères infinitésimaux comme des idéaux maximaux du centre de l'algèbre enveloppante de $\mathrm{Lie} (G)\otimes_\mathbb{R}\mathbb{C}$ ce qui revient au même grâce à l'isomorphisme d'Harish-Chandra.

\

On notera $\pi^{\GL}(\psi)$ la représentation irréductible de $\GL(\mathrm{dim}(\psi),\mathbb{R})$ associé à $\psi$ mais aussi, en fonction du contexte, la représentation tordue comme dans \cite{book} de l'espace tordu associé à $\GL(\mathrm{dim}(\psi),\mathbb{R})$ normalisée à l'aide du modèle de Whittaker comme expliqué en loc. cite. Comme le centralisateur de $\psi$ dans $\GL(n,\mathbb{C})$ est connexe, la représentation $\pi^{\GL}(\psi)$ est l'analogue de $\pi^A(\psi)$ pour le groupe $\GL(\mathrm{dim}(\psi),\mathbb{R})$, ce qui justifie la similitude de la notation.

\

\bf Convention: \rm l'une des difficultés de la théorie est la complexité des notations, on essaie de faire au plus simple, en particulier, ce travail ne concerne que les groupes réductifs réels, on identifiera très souvent $G$ et $G(\mathbb{R})$. Quand on a besoin des points complexes, on revient à la notation $G(\mathbb{C})$. 
\subsection{Endoscopie et A-paquets d'après \cite{book}}\label{defpaquet}
On fixe un paramètre $\psi$ comme dans \ref{parametre}. On note $A(\psi)$ le groupe des  composantes du centralisateur de $\psi$ dans 
$\hat{G}$. On rappelle que l'on a défini la représentation tordue $\pi^{\GL}(\psi)$ en \ref{parametre}. On note $-1_{\mathrm{SL}_2(\mathbb{C})}$ l'élément non trivial du  centre de $\mathrm{SL}_2(\mathbb{C})$ et $s_\psi$ son image par $\psi$.

Soit $\pi$ une combinaison linéaire finie, à coefficients des nombres complexes, de représentations semi-simples irréductibles de $G\times A(\psi)$ et soit $s$ dans le centralisateur dans $\hat{G}$ de $\psi$ tel que $s^2=1$. On note $\pi(s)$ la combinaison linéaire de représentation de $G$ obtenue en évaluant les représentations irré\-ductibles de $A(\psi)$ qui interviennent (ici ce sont des sommes de caractères) en $s$, où  $s$ a été identifié à son image dans $A(\psi)$. On note $H_s$ la donnée endoscopique elliptique de $G$ de la forme $(s,H,\mathcal{H},\xi)$ où $\mathcal{H}$ est le sous-groupe de $^LG$ engendré par la composante neutre du centralisateur de $s$ et par l'image de $\psi$. Ainsi $\psi$ se factorise par $\mathcal{H}$ par construction. On est dans des cas simples qui ne nécessitent pas de données auxiliaires et on voit alors $\psi$ comme un morphisme d'Arthur pour le groupe $H$ de la donnée endoscopique que l'on vient de définir; pour éviter toute ambigüité, on note $\psi_s$ cette factorisation.

 Si $G$ est quasi-déployé, \cite{book} (cf. l'introduction de cette référence) montre qu'il existe une unique combinaison linéaire de représentations unitaires et  semi-simples de $G\times A(\psi)$, notée $\pi^A(\psi)$ telle que 
 
(i) la représentation $\pi^A(\psi)(s_{\psi})$ est stable;

(ii) la trace tordue de $\pi^{\GL}(\psi)$ est le transfert de 
 la trace de $\pi^A(\psi)(s_{\psi})$ pour l'endoscopie tordue;

(iii) pour tout $s$ dans le centralisateur de $\psi$, tel que $s^2=1$, $\pi^A(\psi(s_\psi s))$ est le transfert de la représentation stable de $H_s$, $\pi^A(\psi_s)(s_{\psi_s})$.

En fait Arthur a démontré plus,  si l'on écrit $\pi^A(\psi)$ sous la forme $$
\pi^A(\psi)=:\bigoplus_{\epsilon} \pi^A(\psi,\epsilon)\otimes \epsilon,\eqno(1)$$ où $\epsilon$ parcourt l'ensemble des caractères de $A(\psi)$, on a les propriétés suivantes:

pour tout $\epsilon$, $\pi^A(\psi,\epsilon)$ est une combinaison linéaire à coefficients dans $\mathbb{Z}_{\geq 0}$ de représentations unitaires irréductibles de $G$ et si $\pi^A(\psi,\epsilon)\neq 0$ la restriction de $\epsilon$ au centre de $\hat{G}$ est triviale (ici $G$ est quasi-déployé). 

Cette propriété est équivalente à ce que $\pi^A(\psi)$ soit une représentation semi-simple de longueur finie et unitaire de $G(\mathbb{R})\times A(\psi)$, comme on l'a énoncé dans l'introduction.

A l'opposé si on écrit $\pi^A(\psi)$ sous la forme $\sum_\pi \pi\otimes \epsilon(\pi)$ où $\pi$ parcourt l'ensemble des représentations unitaires irréductibles de $G$, alors \cite{book} ne dit pas grand chose sur les $\epsilon(\pi)$ qui sont non nuls sauf que ce sont des sommes de caractères puisque $A(\psi)$ est commutatif, $A(\psi)$ est même un $2$-groupe. 

Il n'y a aucun contre-exemple à ce que ces $\epsilon(\pi)$ quand ils sont non nuls, soient des caractères (c'est-à-dire des représentations de dimension 1). Pour le moment on ne sait pas démontrer cette propriété pour les groupes réels, alors qu'on sait le faire pour les groupes p-adiques et les groupes complexes. Quand cette propriété  est réalisée on dit que l'on a multiplicité un locale. Cette terminologie est justifiée par le fait que la propriété est équivalente à ce que la restriction de $\pi^A(\psi)$ à $G$ est sans multiplicité.

Ne supposons plus que $G$ est quasi-déployé, alors (iii) ci-dessus permet encore de définir $\pi^A(\psi)$  (cf. \cite{note}); alors que dans le cas quasi-déployé faire $s=1$ dans (iii) ne donne aucune indication. Dans le cas non quasi-déployé $s=1$ est une relation qui complète la définition de $\pi^A(\psi)$.

\subsection{Les cas connus}
Rappelons que l'on connaît $\pi^A(\psi)$ dans le cas où le caractère infinitésimal associé à $\psi$ est entier régulier et  $G$ quasi-déployé. C'est dans \cite{AMR} qui montre que les constructions de \cite{AJ} coïncident avec celles de \cite{book}. 

En \cite{pourhowe} complété en \cite{note} pour le cas non quasi-déployé, on a aussi montré que $\pi^A(\psi)$ a la multiplicité un locale si la restriction de $\psi$ au sous-groupe $\mathbb{C}^*$ de $W_{\mathbb{R}}$ est triviale.

Et on traite le cas très régulier de l'introduction dans \cite{recent} 9.3 sans supposer $G$ quasi-déployé.

\section{Translation et endoscopie}
Contrairement au paragraphe précédent et au but principal de l'article, on démontre la compatibilité de la translation avec l'endoscopie ordinaire et tordue dans un cadre général pour ne plus avoir  à y revenir.
\subsection{Rappel sur le transfert  spectral local\label{transfert}}
On se place dans la situation suivante: ici on remplace $G$ par un espace tordu $\tilde{G}$ et $\bf H$ est une donnée endoscopique elliptique de $\tilde{G}$. On renvoit à \cite{stabilisation} I.1.7 et suivant pour les définitions. On note $\tilde{H}$ l'espace tordu sous-jacent à la donnée endoscopique, la torsion ne peut être qu'une torsion intérieure.

On rappelle que le transfert spectral local, tordu ou non, est défini en toute généralité, en particulier il l'est pour les groupes réels. Le cas tordu, qui inclut le cas non tordu comme cas particulier, se trouve dans \cite{stabilisation} IV.3.3. Bien sûr l'image du transfert est une combinaison linéaire à coefficients complexes d'un nombre fini de représentations irréductibles. On peut difficile\-ment espérer mieux, juste que les coefficients soient des entiers relatifs mais ceci n'est pas démontré.

Ce transfert spectral se traduit par une identité de caractères. Plus précisément soit $\tilde{\pi}^H_{st}$ une combinaison linéaire stable de représentations tordues pour $\tilde{H}$. Alors $\tilde{\pi}$, une combinaison linéaire de représentations tordues pour $\tilde{G}$, est le transfert de $\tilde{\pi}^H_{st}$ si et seulement si pour tout élément semi-simple régulier de $\tilde{G}$, on a l'égalité
\[\tilde{\pi}(g)=\sum_{h\in H/\sim} \Delta(h,g) \tilde{\pi}^H_{st}(h),
\] où les $\Delta(h,g)$ sont les facteurs de transfert (géométriques) définis par Lang\-lands-Shelstad et où on a identifié les représentations à leur caractère. Cette formule est  élémentaire, on renvoit le lecteur à \cite{book} (8.3.4).
\subsection{Translation et endoscopie ordinaire}\label{endoscopie}
Un des outils performants de l'étude des modules de Harish-Chandra pour les groupes de Lie réels est leur mise en familles cohérentes. Pour l'endoscopie, ce n'est pas directement utilisable mais peu s'en faut. 

Dans une situation endoscopique,  la correspondance entre représentations de dimension finie n'est pas le transfert mais simplement l'identité des carac\-tères restreints aux tores des groupes complexifiés et c'est ce qui permet de comprendre le comportement des foncteurs de translation que l'on va définir, vis à vis de l'endoscopie.

Plus précisément, pour tout groupe algébrique $G$ réductif défini sur $\mathbb{R}$, il y a une bijection entre le groupe de Grothendieck des représentations de dimension finie de $G(\mathbb{C})$ et les combinaisons linéaires de caractères algébri\-ques d'un tore maximal de $G(\mathbb{C})$, invariantes sous l'action du groupe de Weyl. 

Soit $\bf H$ une donnée endoscopique elliptique de $G$ dont on note $H$ le groupe sous-jacent. Dans les cas considérés ici, il n'y a pas à faire intervenir de données auxiliaires, toutefois la construc\-tion de ce paragraphe est totalement générale et fonctionne aussi si l'on a des données auxiliaires. Donc $G$ et $H$ sont des groupes algébriques qui sur $\mathbb{C}$ ont des tores maximaux isomorphes, plus précisément les tores duaux sont égaux. Dans le cas où il y a à considérer des données auxiliaires, ici c'est bien le groupe $H$ qui intervient et non son extension centrale. 

On note $T_\mathbb{C}$ un tore algébrique maximal du groupe algébrique complexe $G(\mathbb{C})$, on l'identifie à un tore maximal de  $H(\mathbb{C})$ (cf \cite{stabilisation} 1.1.10). Soit $E$ une combinaison linéaire de caractères algébriques de $T_\mathbb{C}$ que l'on suppose invariante par le groupe de Weyl de $G(\mathbb{C})$; elle est donc aussi invariante pour le groupe de Weyl de $H(\mathbb{C})$. On note alors $E^G$ la combinaison linéaire de représentations de dimension finie de $G$ dont la restriction à $T_\mathbb{C}$ (après prolongement holomorphe à $G(\mathbb{C})$) est $E$ et $E^H$ son analogue pour $H$.

\begin{theorem}
Soit $\pi^H_{st}$ une combinaison linéaire stable de modules de Harish-Chandra pour $H$. On note $\pi_G$ le transfert de $\pi^H_{st}$ à $G$. Soit $E$ comme ci-dessus, alors $\pi^H_{st}\otimes E^H$ est une combinaison linéaire stable de modules de Harish-Chandra pour $H$ dont le transfert à $G$ est $\pi_G\otimes E^G$.
\end{theorem}
On utilise la formule de \ref{transfert}.  Soit $(g,h)$ un couple d'éléments semi-simples et réguliers l'un dans $G$ et l'autre dans $H$. Alors la trace de $g$ dans $E^G$ est égale à la trace de $h$ dans $E^H$ et le théorème s'obtient en multipliant l'égalité de transfert par cette valeur commune.

\subsection{Translation et endoscopie tordue }\label{endtord}
Attention, ici on change les notations par rapport à \ref{transfert} et donc par rapport au paragraphe précédent, on n'a pas trop le choix, il faut à un moment inverser les notations quand on travaille avec un groupe, $G$, donné. On regarde les représentations de $G$ qui sont obtenues par transfert via ses données endoscopiques elliptiques et on regarde le transfert des combinaisons linéaires stables de représentations de $G$ dans les situation où $G$ est le groupe d'une donnée endoscopique elliptique  tordue. 

 On fixe donc un espace tordu $\tilde{H}$, on peut même introduire un caractère $\omega$ comme dans \cite{stabilisation}, on note $H$ le groupe sous-jacent à $\tilde{H}$. On reprend l'article de Bergeron et Clozel \cite{BC} pour y mettre l'endoscopie.

 Ici $G$ est le groupe sous-jacent à une donnée endoscopique de l'espace tordu. On a encore une identification, sur le corps des complexes, entre un tore maximal $T_{G,\mathbb{C}}$ de $G(\mathbb{C})$ et $\tilde{T}_{H,\mathbb{C}}/(1-\theta)T_{H,\mathbb{C}}$ où $\tilde{T}_{H,\mathbb{C}}$ est un tore tordu maximal  de $\tilde{H}(\mathbb{C})$ dont $T_{H,\mathbb{C}}$ est le groupe sous-jacent. Ici $\theta$ est un élément de $\tilde{H}$ convenable, fixant un épinglage de $H(\mathbb{C})$.  On a alors les morphismes:
$$
\tilde{T}_{H,\mathbb{C}} \rightarrow \tilde{T}_{H,\mathbb{C}}/(1-\theta)T_{H,\mathbb{C}} \simeq \tilde{T}_{G,\mathbb{C}};\eqno(1)$$
 ici $\tilde{T}_{G,\mathbb{C}}$ est simplement le tore complexe $T_{G,\mathbb{C}}$ si la donnée endoscopique $G$ n'est pas à torsion intérieure propre et sinon c'est une trivialisation d'un tore tordu de cette donnée. On n'insiste pas (en renvoyant le lecteur au premier chapitre de \cite{stabilisation}) car on ne s'intéresse qu'au cas où la donnée endoscopique n'est pas tordue et on oubliera le $\tilde{T}_{G,\mathbb{C}}$ pour le remplacer par $T_{G,\mathbb{C}}
$. Notons $N$ le composé de ces morphismes et la propriété importante est que si $\tilde{h}\in \tilde{H}$ et $g\in G$ sont des éléments semi-simples réguliers tels que le facteur de transfert géométrique $\Delta(g,\tilde{h})$ est non nul alors il existe $h_{\mathbb{C}}\in \tilde{T}_{H,\mathbb{C}}$ tel que $\tilde{h}$ est conjugué sous $H(\mathbb{C})$ de $\tilde{h}_{\mathbb{C}}$ et $g$ est conjugué sous $G(\mathbb{C})$ de $N(\tilde{h}_\mathbb{C})$.

On veut relier le groupe de Grothendieck des représentations tordues et de dimension finie de $\tilde{H}$ et le groupe de Grothendieck des représentations de dimension finie de $G$. Comme en \ref{endoscopie}, on passe par les restrictions aux tores et il faut donc contrôler l'action des groupes de Weyl.

Il y a un cas très joli qui permet d'établir une bijection entre ces groupes de Grothendieck. C'est le cas où $G$ est le groupe sous-jacent d'une donnée endoscopique principale. 

On traite donc d'abord ce cas, c'est à dire que le $s$ de la donnée endoscopi\-que est $\theta^*$, l'automorphisme du groupe dual de $H$ dual de $\theta$. 

\begin{theorem} (avec les notations précédentes) On suppose que $G$ est une donnée endoscopique principale (non tordue) de l'espace tordu $\tilde{H}$. 

 Soit $E$ une combinaison linéaire $W_G$ invariante de caractères algébriques d'un tore maximal de $G$. Alors il existe une combinaison linéaire de représen\-tations tordues de dimension finie, $E^{\tilde{H}}$, de $\tilde{H}$ et une combinaison linéaire de représentations de dimension finie $E^G$ de $G$, telles que pour  tout élément semi-simple régulier $\tilde{h}\in \tilde{T}_H$ on ait
$$
\mathrm{trace}_{E^{\tilde{H}}}(\tilde{h})= \mathrm{trace}_E(N(\tilde{h}))=\mathrm{trace}_{E^G}(N(\tilde{h})).\eqno(2)
$$

\end{theorem}
Pour démontrer le théorème, il suffit de le faire dans le cas élémentaire où $E$ est obtenu en sommant sous l'action de $W_G$ un caractère algébrique dominant $\nu$, c'est-à-dire que $E=\sum_{w\in W_G}w.\nu$.

On note $\mu$ le caractère de $T_H$ qui s'obtient en relevant $\nu$; c'est un caractère algébrique et dominant.

On note $E_\mu$ la représentation de $H$ de plus haut poids $\mu$. On en fait une représentation tordue de $\tilde{H}(\mathbb{C})$ en imposant que $\theta$ opère trivialement sur l'espace de poids $\mu$. Ainsi, le tore tordu $\tilde{T}_H$ opère sur la somme des vecteurs de poids extrémaux, représentation notée $\tilde{E}_{ext}$ et il existe une représentation tordue $E^{\tilde{H}}$ de $\tilde{H}$ combinaison linéaire de représentations de dimension finie dont la restriction au tore tordu est précisément $\tilde{E}_{ext}$. On va montrer que $E^{\tilde{H}}$ répond à la condition du théorème.

Pour tout $w\in W_H$, le groupe de Weyl de $H$, on fixe $n_w\in H$ un représentant de $w$. Soit $v$ un vecteur non nul de $E_\mu$ de poids $\mu$ et $\tilde{E}_{ext}$ se réalise dans l'espace vectoriel
$$
\oplus_{w}\mathbb{C}\, n_w .v, \eqno(3)
$$
où $w$ parcourt un système de représentants de $W_H/W_{H,\mu}$ où $W_{H, \mu}$ est le stabilisateur de $\mu$ dans $W_H$. On remarque d'abord que $\mathbb{C}\, n_w.v$ est un espace invariant par $\theta$ si et seulement si $w\in (W_H/W_{H,\mu})^\theta$. Comme $\mu$ est dominant, $W_{H,\mu}$ est le groupe de Weyl d'un sous-groupe de Levi de $H$ et comme $\mu$ est $\theta$ invariant, ce sous-groupe de Levi est stable sous $\theta$. Ainsi le représentant de Kostant de $w$ dans la classe 
$(W_H/W_{H,\mu})$ est nécessairement invariant sous $\theta$.

On doit donc calculer l'action de $\theta$ sur tout espace vectoriel $\mathbb{C}n_w.v$ en supposant que $w\in W_H^{\theta}$. On sait que $W_H^\theta$ s'identifie à $W_G$, c'est l'hypothèse que $G$ fait partie de la donnée endoscopique principale de $\tilde{H}$ et cela de façon compatible au morphisme (1). Mais on sait aussi que $W_H^\theta$ s'identifie au groupe de Weyl de $H(\mathbb{C})^{\theta,0}$ d'après \cite{KS} avant dernier paragraphe avant (1.2). On fixe donc $n^0_w$ dans $H(\mathbb{C})^\theta$ représentant $w$. Or pour $w$ et $n_w$ comme ci-dessus, on a
$$
\theta (n_w).v=n_w. (n_w^{-1}.\theta(n_w)).v=n_w v
$$D'où $n_w^{-1}.\theta(n_w)=t\theta(t)^{-1}$ où $t=n_w^{-1}n_w^0$.
Ainsi $\theta$ agit trivialement  $n_w.v$ et par compatibilité des actions de $W_H^\theta$ et $W_G$ l'action de $\tilde{T}_H$ sur cet espace est l'image réciproque de $w.\nu$ par (1). C'est exactement ce qu'il fallait démontrer pour obtenir (2). D'où le théorème.

\

En fait on a aussi besoin du cas où la donnée endoscopique n'est pas principale et ce qui change, alors, est que le groupe de Weyl de la donnée endoscopique s'identifie à un sous-groupe de $W_H^ \theta$.

On a alors un théorème analogue à celui de \ref{endoscopie}.
\begin{theorem}

 Soit $E$ une combinaison linéaire $W_H ^\theta$ invariante de caractères algébriques d'un tore tordu maximal de $\tilde{H}$. Alors il existe une combinaison linéaire de représentations tordues de dimension finie, $E^{\tilde{H}}$, de $\tilde{H}$ et une combinaison linéaire de représentations de dimension finie $E^G$ de $G$, telles que pour  tout élément semi-simple régulier $\tilde{h}\in \tilde{T}_H$ on ait
$$
\mathrm{trace}_{E^{\tilde{H}}}(\tilde{h})= \mathrm{trace}_{E^G}(N(\tilde{h}))=\mathrm{trace}_E(N(\tilde{h})s.\eqno(2)
$$

\end{theorem}
C'est la m\^eme démonstration que ci-dessus mais en partant d'un poids dominant $\mu$ pour $H(\mathbb{C})$, $\theta$-invariant.

\

\begin {theorem} Soit $E$ comme dans le théorème précédent dont on reprend les hypothèses et notations. Et soit $\pi_G$ une combinaison linéaire stable de modules de Harish-Chandra de $G$ (ou pour une donnée auxiliaire de la donnée endoscopique). Alors on a l'égalité de transfert:
\[\mathrm{transfert}\left(\pi_G\otimes E^G \right)=\left(\mathrm{transfert}(\pi_G)\right)\otimes E^{\tilde{H}}.
\]

\end{theorem} 
La démonstration est exactement la même que celle de \ref{endoscopie} en tenant compte de ce qui précède et de \ref{transfert}.
\subsection{Caractère infinitésimal, transfert et translation}\label{translation}
Un caractère infinitésimal sera maintenant (cf. introduction) la donnée d'un idéal maximal du centre de l'algèbre enveloppante, $Z(\mathfrak{g})$ de $G$. Une telle donnée est uniquement déterminée par une orbite du groupe de Weyl de $G(\mathbb{C})$ dans l'espace vectoriel dual de $\mathfrak{T}$, l'algèbre de Lie complexifiée d'un tore maximal de $G$. En reprenant aussi la donnée endoscopique $\bf H$ du paragraphe \ref{endoscopie}, on voit qu'il existe une application de l'ensemble des caractères infinitésimaux pour $H$ dans l'ensemble des caractères infinitésimaux pour $G$ et cette application n'est en général pas injective.

Soit $V_G$ un module de Harish-Chandra de longueur finie pour $G$. Il existe une décomposition en somme directe finie
\[V_G=\oplus_\nu V_{G,\nu}\]
où $\nu$ parcourt l'ensemble des caractères infinitésimaux pour $G$ et où $V_{G,\nu}$ est l'ensemble des éléments de $V_G$ annulés par une puissance de l'idéal de $Z(\mathfrak{g})$ défini par $\nu$. Reprenons la donnée endoscopique $\bf H$ et soit $\nu$ un caractère infinitésimal pour $G$. Soit $V_H$ un module de Harish-Chandra pour $H$, on note $V_{H,\nu}$ la somme des $V_{H,\nu_H}$ où $\nu_H$ parcourt l'ensemble des caractères infinitésimaux pour $H$ se transférant en $\nu$; bien sûr le transfert d'un caractère infinitésimal vient de l'inclusion des $L$-groupes. On prolonge  ces définitions par linéarité pour les transporter aux combinaisons linéaires de modules de Harish-Chandra. On a alors la propriété de base de l'endoscopie:
\begin{rmq}Soit $V_H$ une combinaison stable de modules de Harish-Chandra pour $H$ et soit $V_G$ son transfert à $G$. Pour tout caractère infinitésimal $\nu$ pour $G$, $V_{H,\nu}$ est stable et a pour transfert $V_{G,\nu}$.
\end{rmq}

De m\^eme si $G$ est une donnée endoscopique  d'un espace tordu $\tilde{H}$, on a une application de l'ensemble des caractères infinitésimaux pour $G$ dans l'ensemble des caractères infinitésimaux pour $H$ qui n'est pas injective en général. 

On appliquera cette remarque en la couplant avec  \ref{endoscopie} (et \ref{endtord} pour le cas tordu) pour obtenir avec les notations de \ref{endoscopie}:
\[\mathrm{transfert}\biggl((V_H\otimes E^H)_\nu\biggr)=(V_G\otimes E^G)_\nu.\]

\subsection{Translation et caractère infinitésimal, une remarque}\label{infinitesimal}
On rappelle un résultat connu mais non trivial de \cite{HS}. Soient $\pi$ un module de Harish-Chandra pour $G$ ayant un caractère infinitésimal, noté $\nu$. Et soit  $E$ une combinaison linéaire de représentations de dimension finie de $G$ et soit $\mathcal{E}$ l'ensemble des caractères algébriques d'un tore maximal de $G(\mathbb{C})$ opérant dans $E$. L'ensemble $\nu+\mu$ où $\mu$ parcourt $\mathcal{E}$ est bien défini comme un ensemble de caractères infinitésimaux. Alors les modules de Harish-Chandra irréductibles intervenant dans la combinaison linéaire $\pi\otimes E$ ont tous leur caractère infinitésimal dans l'ensemble que l'on vient de définir.

\section{Applications aux paquets d'Arthur\label{translationpaquet}}
\subsection{Première forme d'utilisation de la translation\label{premiereforme}}
On revient à la situation des groupes classiques.
On fixe un paramètre $\psi$ comme en \ref{parametre};  on suppose que $\psi$ est de bonne parité. On décompose encore $\psi$ en $$\bigoplus_{i\in[1,v]}\rho_i\otimes R[a_i] \bigoplus \psi_u$$ où les $\rho_i$ sont des séries discrètes de paramètre respectif $t_i$, les $a_i$ sont des entiers et où $\psi_u$ est la partie unipotente de $\psi$. On suppose que $t_1 \geq \cdots \geq t_v$.
On fixe un autre paramètre $\psi_+$ et on dit (comme dans l'introduction) que $\psi_+$ est très dominant par rapport à $\psi$ si  
\[\psi_+=\bigoplus_{i\in[1,v]}\rho'_i\otimes R[a_i] \bigoplus \psi_u\]
donc la seule chose qui change par rapport à $\psi$ sont les séries discrètes $\rho'_i$ dont le paramètre est noté $t'_i$ et qui doit vérifier \[T_1:=t'_1-t_1 >> \cdots >> T_v:=t'_v-t_v>>0.\]
On rappelle aussi qu'à $\psi$ est associé un caractère infinitésimal que l'on note $\nu_\psi$.

On réalise la forme quasi-déployée de $G$ comme le groupe sous-jacent d'une donnée endoscopique elliptique d'un espace tordu $\tilde{H}$ où $H$ est $\GL(n^*)$ avec $n^* $ la dimension de la représentation $\psi$.

On considère la représentation de dimension finie de $\GL(n^*)$ de plus bas  poids $-\lambda$ où $$\lambda:=T_1, \cdots, T_1, T_2, \cdots,T_2, \cdots, T_v, \cdots, T_v, 0, \cdots,0,-T_v,\cdots, -T_1)$$ où chaque $T_i$ et chaque $-T_i$ pour $i\in [1,v]$ vient $a_i$ fois. Et on note $E(\psi_+,\psi)$ la représentation tordue pour le tore diagonale tordu de $\tilde{\GL}(n^*)$ agissant dans la somme des espaces de poids extr\^emaux. On a donc défini les combinaisons linéaires de représentations de dimension finie, $E^{\GL}(\psi_+,\psi)$ (on oublie le $\tilde{\empty}$ dans la notation), $E^G(\psi_+,\psi)$. On voit aussi $\lambda$ comme un poids pour $G$ en ne prenant que les $r$ premières coordonnées où $r$ est le rang déployé de $G(\mathbb{C})$ et on note simplement $E^G_+$ la combinaison linéaire de représentations de dimension finie de $G$ dont la restriction à un tore maximal de $G(\mathbb{C})$ est la somme des caractères dans l'orbite de $-\lambda$ sous l'orbite du groupe de Weyl de $G$.

\

On simplifie les notations en notant $E_+$ à la place de $E(\psi_+,\psi)$ que l'on voit aussi comme une  combinaison linéaire des caractères d'un tore maximal de $G(\mathbb{C})$.  Et on note $E^{G}_+$ au lieu de $E^G(\psi_+,\psi)$.

Les groupes orthogonaux impairs et symplectiques quasi-déployés sont des données endoscopiques principales de l'espace tordu et on a décrit $E^G_+$ dans \ref{endtord}. Si $G$ est un groupe spécial orthogonal pair,  le groupe de Weyl de ce groupe est plus petit que le groupe de Weyl de l'espace tordu. Mais  on travaille en fait avec à conjugaison près par l'automorphisme extérieur qui vient du groupe orthogonal, cela est normal d'ajouter cet automorphisme au groupe de Weyl. Ainsi on travaille avec $E^G_+$ plutôt que simplement avec l'orbite sous le groupe de Weyl du groupe spécial orthogonal de $-\lambda$.

\begin{theorem} On a la formule de translation:
\[\pi^A(\psi)= \biggl(\pi^A(\psi_+)\otimes E^G_+\biggr)_{\nu_\psi}\eqno(1)\]
\end{theorem}
 
On commence par le cas où $G$ est quasi-déployé. Ce cas est  plus difficile, car il y a plus de propriétés à vérifier à cause de l'induction tordue.

 Pour démontrer (1), il faut montrer que le membre de droite de (1) vérifie les trois conditions de \ref{defpaquet}. Montrons (i): la stabilité se voit par le fait que le caractère est invariant par conjugaison sous $G(\mathbb{C})$. Cette invariance est vraie pour les représentations de dimension finie donc elle l'est pour le produit tensoriel avec $\pi^A(\psi_+)(s_\psi)$ puisque cela est vrai pour $\pi^A(\psi_+)(s_\psi)$.  Il reste à voir que la projection sur les espaces propres généralisés sous l'action du centre de l'algèbre enveloppante préserve la stabilité. On utilise une autre caractérisation de la stabilité. Une distribution est stable si elle est annulée par les fonctions dont le transfert stable au groupe lui-même est nul. Or ce transfert est compatible à l'action du centre de l'algèbre enveloppante (cf. \cite{stabilisation} premier lemme de la page 443 où l'on fait $\tilde{G}=G=G'$). Cela termine la preuve de (i).

Montrons (ii). D'après \ref{translation}, le transfert tordu du membre de droite de (1)  est 
\[\biggl(\pi^{{\GL}}(\psi_+)\otimes E^{\GL}_+\biggr)_{\nu_\psi}. \eqno(2)
\]
On va  identifier (2) à $\pi^{\GL}(\psi)$. On note encore $n^*$ la dimension de la représentation $\psi$.  On a déjà défini $\lambda$ avant l'énoncé comme  caractère du tore de $\GL(n^*,\mathbb{C})$ (ou plutôt sa différentielle); on note $\lambda'$ la différentielle de $\lambda$. Les caractères infinitésimaux intervenant dans $\pi^{{\GL}}(\psi_+)\otimes E^{\GL}_+$ sont de la forme
\[\nu_{\psi_+}+w\lambda',\]
où $w$ parcourt le groupe des symétries, $\theta$-invariantes, de $n^*$ éléments.
En prenant des représentants des classes de conjugaison sous l'action de ce groupe de symétries, on peut écrire $\nu_{\psi_+}$
sous la forme d'une collection de nombres ordonnés de la forme
$$\nu_+:=
t'_1+(a_1-1)/2, \cdots, t'_1-(a_1-1)/2, \cdots, t'_v-(a_v-1)/2,$$
$$ \nu_{\psi_u},$$
$$ -t'_v+(a_v-1)/2, \cdots, 
-t'_1+(a_1-1)/2, \cdots, -t'_1-(a_1-1)/2.
$$
On veut transformer cette collection de nombres en son analogue où les $t'_i$ sont remplacés par les $t_i$ pour tout $i\in [1,v]$ en utilisant les coordonnées de $\lambda'$. On remarque d'abord  que si l'on n'enlève pas $T_1$ à chacune des $a_1$ premières coordonnées, on enlèvera à au moins une de ces coordonnées un $T_i$ ou $-T_i$ (avec $i\in [1,v]$) ou $0$ et dans tous les cas on trouvera in fine au moins une coordonnée qui n'apparaît pas dans $\nu_\psi$. De même il faut ajouter $T_1$ aux $a_1$ dernières coordonnées. Puis de proche en proche on vérifie que seul $\nu_{+}-\lambda'$ est dans l'orbite de $\nu_\psi$.

Ensuite on utilise le fait que $\pi^{\GL}(\psi)$ est l'induite cohomologique du produit d'un caractère   de $\times_{i\in [1,v]}\GL(a_i,\mathbb{C})$ avec la représentation $\pi^A(\psi_u)$ du groupe $\GL(n^*_u,\mathbb{R})$ où $n^*_u $ est la dimension de $\psi_u$. Le caractère des premiers facteurs vaut $\times_i(det/{\overline{det}})^{t'_i}$. L'induction cohomologique  commute à la tensorisation par une représentation de dimension finie.  Et d'après ce que l'on vient de voir (2) est donc l'induite cohomologique de même type où l'on change simplement le caractère des premiers facteurs en remplaçant $t'_i$ par $t_i$ pour tout $i\in [1,v]$. Et on obtient $\pi^{\GL}(\psi)$. Si $G$ fait partie de la donnée endoscopique principale, on a terminé puisque dans \ref{translation} la projection sur le caractère infinitésimal du côté endosco\-pique se fait uniquement sur $\nu_\psi$. Sinon, $G$ est un groupe spécial orthogonal pair et il y a deux caractères infinitésimaux à considérer qui se déduisent l'un de l'autre par l'action du groupe orthogonal. Mais dans la définition même de $\pi^A(\psi)$ on a considéré le groupe orthogonal et non le groupe spécial orthogonal.

\

Montrons  la propriété (iii) de \ref{defpaquet} pour le membre de droite de (1). On fixe une donnée endoscopique elliptique, d'où un élément semi-simple $s$ du groupe dual de $G$. On suppose que $s^2=1$ et  que $s$ se trouve  dans le centralisateur de $\psi_+$. Ainsi $s$ est aussi dans le centralisateur de $\psi$ et on obtient par cette construction tous les éléments du centralisateur de $\psi$ d'ordre deux. On a aussi $H$ le groupe de la donnée. En partant de $E_+$, on a défini $E_+^H$. On a aussi $\psi_s$ et $\psi_{+,s}$. On a par définition les égalités de transfert qui réalisent $\pi^A(\psi)(s_\psi s)$ comme transfert de $\pi^A(\psi_s)(s_{\psi_s})$ et l'analogue en remplaçant $\psi$ par $\psi_+$. En tenant compte de \ref{translation}, il faut démontrer que
$$\left(\pi^A(\psi_{s,+})(s_{\psi_s})\otimes E_+^H\right)_{\nu_\psi}= \pi^A(\psi_s)(s_{\psi_s}).
  $$
Pour démontrer cette égalité, il suffit encore de la démontrer après transfert à l'espace tordu dont $H$ est une donnée endoscopique; cet espace tordu est un produit de deux espaces tordus de groupe sous-jacent $\GL(n')\times \GL(n'')$ avec $n'+n''=n^*$. En effet si l'égalité est vraie après transfert, elle est vraie pour les distributions stables que l'on a transférées car le transfert endoscopique tordu des fonctions lisses $K$-finies sur  $\tilde{\GL}(m,\mathbb{R})$ vers les fonctions lisses $K$-finie d'un groupe classique qui fait partie d'une donnée endoscopique elliptique de $\tilde{GL}(m,\mathbb{R})$ (pour tout entier $m$) est surjective, modulo bien évidemment l'espace des fonctions dont les intégrales orbitales stables sont nulles et comme toujours modulo l'action du groupe orthogonal pour les groupes orthogonaux pairs. Une référence totalement générale pour cela est la proposition de I.4.11 (ii) et IV.3.4 (iii) (pour la préservation de la $K$-finitude) de \cite{stabilisation}; d'après cette référence qui caractérise l'image du transfert tordu en toute généralité, il faut montrer que l'action des automorphismes de la donnée endoscopique est trivial sur l'espace d'arrivée. Ce groupe d'auto\-morphismes est décrit en \cite {stabilisation} I.1.5 comme une extension des automorphi\-smes  \og extérieurs\fg{} par le groupe $C^{Gal(\mathbb{C}/\mathbb{R})}$ où  $C:=(Z(\hat{H})/Z(\hat{H})\cap Z(\hat{G}))$. Le groupe des automorphismes \og extérieurs\fg{} est trivial sauf dans le cas des groupes orthogonaux pairs où il s'identife au quotient du groupe orthogonal par le groupe spécial orthogonal. Son action sur les fonctions est l'action naturelle. Comme évidemment $Gal(\mathbb{C}/\mathbb{R})$ agit trivialement sur $C$, on doit calculer l'action de $C$. Mais celle-ci est triviale (cf. bas de la page 40 de \cite{stabilisation} où l'action est décrite). D'où l'assertion de surjectivité.   

 On se ramène donc à un espace tordu dont le groupe est un produit de deux groupes linéaires généraux, comme ci-dessus. Le transfert de $\nu_\psi$ à cet espace tordu n'est pas une seule orbite sous le groupe de Weyl de l'espace tordu, en fait c'est une seule orbite sous le groupe de Weyl de $\GL(n^*)$ quand on voit $\GL(n')\times \GL(n'')$ comme un sous-groupe de Levi de $\GL(n^*)$. La même démonstration que celle que l'on a donnée ci-dessus donne le résultat cherché.

 Et cela termine la preuve dans le cas quasi-déployé. Dans le cas général, il n'y a que (iii) à montrer et cela se fait aisément comme ci-dessus. 
\subsection{Remarque sur le théorème précédent\label{translationrep}}

On reprend les notations $\pi^A(\psi,\epsilon)$ de \ref{defpaquet}. 
\begin {cor} On fixe un paramètre $\psi$ de bonne parité et on choisit $\psi_+$ très dominant par rapport à $\psi$. Pour tout caractère $\epsilon$ du centralisateur de $\psi$ dans $\hat{G}$, on a 
\[\pi^A(\psi)(\epsilon)=\biggl(\pi^A(\psi_+)(\epsilon)\otimes E^G(\psi_+,\psi)\biggr)_{\nu_\psi}\eqno(1)
\]
\end {cor}
Comme dans la preuve de la propriété (iii) en \ref{premiereforme} on identifie les éléments du centralisateur de $\psi_+$ d'ordre deux avec leurs analogues pour le centralisa\-teur de $\psi$. Et comme on l'a vu dans la preuve, l'égalité du théorème de \ref{premiereforme} est une égalité de représentations de $G\times A(\psi_+)$, $A(\psi)$ étant vu comme un quotient de $A(\psi_+)$.
Cette application de quotient est non  injective en général. En particulier, on peut ajouter au corollaire que si $\epsilon_+$ est un caractère du groupe des composantes du centralisa\-teur de $\psi_+$ non trivial sur le noyau de l'application de $A(\psi_+)$ sur $A(\psi)$ alors le membre de droite
de (1) (avec $\epsilon_+$) est nul. Montrons ces assertions.
On voit $\pi^A(\psi)$ et $\pi^A(\psi_+)$ comme des représentations de $G$ fois le centralisateur de $\psi_+$. Ainsi on a montré dans le théorème l'égalité de $\pi^A(\psi)$ avec  
$$
\biggl(\pi^A(\psi_+)\otimes E^G(\psi_+,\psi)\biggr)_{\nu_\psi}. \eqno(2)
$$
Comme $\pi^A(\psi)$ ne fait intervenir que les représentations du centralisateur de $\psi_+$ qui ne se factorisent pas par la surjection de $A(\psi_+)$ sur $A(\psi)$. Il en est de même du terme (2). Et on obtient alors facilement les assertions cherchées.

\subsection{Calcul de la translation, préliminaire}\label{calcul}
Evidemment en général on ne sait pas calculer la translation de \ref{translation}. Toutefois la translation utilisée est très bien adaptée aux représentations obtenues par induction cohomologique.  On fixe $L,\mathfrak{q}$ un couple formé d'un sous-groupe de $G$ défini sur $\mathbb{R}$ tel que $L(\mathbb{C})$ soit un sous-groupe de Levi de $G(\mathbb{C})$ pour l'algèbre parabolique $\mathfrak{q}$. On remarque que $L$ est un produit de  groupes unitaires (toutes les formes interviennent) par un groupe $G_0$, spécial orthogonal ou symplectique non nécessairement quasi-déployé.

On reprend les notations $v, t_i,t'_i$ de \ref{translation}. Et on suppose qu'il y a $v$ groupes unitaires dans le produit, c'est-à-dire que $t_i\neq 0$ pour $i\in [1,v]$ et $t_i=0$ si $i>v$.

On note $\chi$ le caractère de ce produit de groupes unitaires qui moralement correspond aux $t_i$ pour $i\in [1,v]$; il y a une subtilité puisque si l'un des $t_i$ est demi-entier, il ne lui correspond de caractère que pour un revêtement d'un groupe unitaire. Mais comme on veut une induction cohomologique qui préserve le caractère infinitésimal, on doit aussi prendre une racine carré de la  demi-somme des racines du radical nilpotent de $\mathfrak{q}$. Il n'est pas miraculeux que les deux difficultés s'annulent; on induit cohomologiquement le caractère qui sur $U(p_i,q_i)$ vaut $det^{\tilde{t}_i}$ où $$\tilde{t}_i=t_i+(a_i-1)/2+\epsilon_G+\sum_{j>i}a_j+n_0,
$$
où $n_0$ est le rang de $G_0(\mathbb{C})$ et où $\epsilon_G=0$ si $G$ est un groupe spécial orthogonal pair, $1/2$ si $G$ est un groupe spécial orthogonal impair et $1$ si $G$ est un groupe symplectique. Et on vérifie que l'assertion de bonne parité est équivalent à ce que $t_i+(a_i-1)/2$ est un entier si $G$ n'est pas un groupe orthogonal impair et un demi-entier sinon. Donc dans tous les cas $\tilde{t}_i$ est un entier et permet donc de définir un caractère.

On définit de façon analogue $\chi_+$ en utilisant les $t'_i$.

On fixe $\sigma$ un module de Harish-Chandra, irréductible et unitaire, pour $G_0$. On peut donc former les induites cohomologiques $
A_{\mathfrak{q}}(\chi_+\times \sigma)$ et $A_{\mathfrak{q}}(\chi \times \sigma)$; ici on n'est pas très précis sur la définition (il y a suffisamment de références dans la littérature par exemple \cite{KV}) et on renvoit le lecteur à \cite{recent}. Il nous suffit, ici, de dire que l'on suppose que $\sigma$ est faiblement unipotent (cf. ci-dessous \ref{unipotent}) de caractère infinitésimal le même que celui associé à $\psi_u$ et que l'on utilise la définition pour l'induction cohomologique telle que le caractère infinitésimal de l'induite cohomologique est le caractère infinitésimal associé à $\psi$; comme les groupes unitaires ont des tores connexes, cela fixe la situation à la torsion près par un caractère du facteur qui est un groupe orthogonal (quand il y en a un) et cela n'importe pas pour ce qui suit.

L'induction cohomologique est a priori un foncteur dérivé qui vit donc en plusieurs degrés. On va voir que dans notre cas, il est nul en tout degré sauf au plus un, ce qui justifie la notation $A_{\mathfrak{q}}()$:

 on a supposé  que $\sigma$ est faiblement unipotent  et on suppose en plus   que $t_1\geq \cdots \geq t_v$. Comme $\chi_+$ est supposé très grand par rapport au caractère infinitésimal de $\sigma$ la représentation $
A_{\mathfrak{q}}(\chi_+\times \sigma)$ se trouve dans le bon domaine, c'est-à-dire le good range (cf. \cite{KV}), et cette représentation est donc unitaire et irréductible. Par contre la représentation $
A_{\mathfrak{q}}(\chi\times \sigma)$ est seulement dans le weakly fair range et cette représentation est donc unitaire mais pas nécessairement irréductible (cf. le chapitre 12 de \cite{KV} et en particulier le dernier théorème 12.9). 

On reprend la notation $E^G_+$ qui simplifie $E^G(\psi_+,\psi)$ de \ref{translation}.
\begin{theorem} On a:
 $$A_{\mathfrak{q}}(\chi\times \sigma)=\left( A_{\mathfrak{q}}(\chi_+\times \sigma)\otimes E^G_+ \right)_{\nu_\psi}.
 $$
\end{theorem}
Pour démontrer ce théorème on écrit qu'après semi-simplification 
$A_{\mathfrak{q}}(\chi_+\times \sigma)\otimes E^G_+$ est isomorphe à $A_{\mathfrak{q}}(\rho)$ où $\rho$ est le produit tensoriel de $\chi_+ \times \sigma$ avec la restriction de $E^G_+$ à $L$. On écrit cette restriction. On sait que la restriction de $E^G_+$ à un tore de $L$ est $\sum_w w.\lambda^{-1}_+$ (en écriture multiplicative) où $\lambda_+$ est le caractère décrit en \ref{translation} et où $w$ parcourt le groupe de Weyl, $W^G$,  de $G$. On décompose suivant les classes de $W^G$ modulo le groupe de Weyl de $L$ opérant à gauche (avec nos notations). Chaque fois que l'on fixe une telle classe, on obtient une combinaison linéaire de représentations de dimension finie pour $L$ dont la restriction à un tore de $L$ est $\sum_{w\in W_L}ww_0 \lambda^{-1}_+$ où $w_0$ est un représentant de cette classe. Les seuls caractères infinitésimaux qui après induction de $L$ à $G$ sont conjugués de $\nu_\psi$  sont alors pour la classe de l'identité de $W^G$. D'où le théorème.

\subsection{Le cas unipotent}\label{unipotent}
En suivant Barbasch et Vogan, on a appelé morphisme spécial unipotent tout morphisme $\psi$ de bonne  parité dont la restriction au sous-groupe $\mathbb{C}^*$ de $W_\mathbb{R}$ est trivial, c'est-à-dire qui se factorise par le groupe de Galois de $\mathbb{C}/\mathbb{R}$. Dans la théorie de \cite{ABV}, les représentations associées à un tel morphisme sont  faiblement unipotentes au sens suivant que l'on tire de \cite{KV} 12.3 (on utilise ici le fait que les groupes que l'on considère sont semi-simples):

soit $\pi$ un module de Harish-Chandra, irréductible, pour $G$ dont on note $\nu_\pi$ le caractère infinitésimal. On dit que $\pi$ est faiblement unipotent si $\nu_\pi$ est réel (c'est-à-dire) dans le $\mathbb{R}$-espace vectoriel engendré par l'ensemble des racines de $G$ et si  pour toute représentation de dimension finie $E$ de $G$ et pour tout caractère infinitésimal $\nu$, la représentation
$\biggl(\pi\otimes E\biggr)_\nu$ est non nulle seulement si $\vert \nu \vert$ est supérieur ou égal à $\vert \nu_\pi\vert$. 
\begin{theorem}
Soit $\psi$ un paramètre unipotent, alors $\pi^A(\psi)$  est une combinaison linéaire de représentations faiblement unipotentes de $G$ à coefficients dans l'espace des caractères de $A(\psi)$.
\end{theorem}
La condition sur le caractère infinitésimal qui doit être réel est clairement vérifiée.

Ensuite il suffit de démontrer que la représentation $\pi^A(\psi)(s_\psi)$ est faible\-ment unipotente: en effet l'image de $s_\psi$ dans le groupe des composantes du centralisateur de $\psi$ coïncide avec l'image d'un élément central. Or les caractères du groupe des composantes intervenant dans la description de $\pi^A(\psi)$ prennent la même valeur sur tout élément du centre de $\hat{G}$. Comme les représentations $\pi^A(\epsilon)$ de \ref{defpaquet} sont des combinaison linéaires à coefficients positifs ou nuls de représentations irréductibles (cf. \cite{note} pour le cas non quasi-déployé qui généralise le cas quasi-déployé dû à Arthur), en évaluant en $s_\psi$, il n'y a pas de simplifications possibles. D'où la réduction annoncée. 

Pour démontrer la propriété de minimalité du caractère infinitésimal de $\pi^A(\psi)$ pour les tensorisations $\pi^A(\psi)\otimes E$, il suffit évidemment de faire parcourir à $E$ l'ensemble des combinaisons linéaires de représentations de dimension finie. On peut donc plutôt supposer que $E$ est une combinaison linéaire de caractères algébriques d'un tore de $G(\mathbb{C})$ invariante sous le groupe de Weyl et noter $E^G$ la combinaison linéaire de représentations de dimension finie de $G$ ayant $E$ pour restriction au tore. On note alors $E^{\GL}$ la représenta\-tion de l'espace tordu $\tilde{\GL}(n^*,\mathbb{R})$ (ou $n^*=\mathrm{dim}\, \psi$) qui correspond elle aussi à $E$ (cf. \ref{endtord}). On a alors vu en \ref{translation} que pour tout caractère infinitésimal $\nu$ pour l'espace tordu le transfert de $$\mathrm{transfert}
\left(  \pi^A(\psi)(s_\psi)\otimes E^G \right)_\nu =
\left( \pi^{\GL}(\psi) \otimes E^{\GL}\right)_\nu.
$$
Le terme de gauche est nul si celui de droite l'est parce que le transfert géométrique est surjectif dans l'ensemble des fonctions \og lisses\fg{} modulo l'ensemble des fonctions dont les intégrales orbitales stables sont nulles (on a donné la preuve de cette surjectivité en \ref{premiereforme}). Or du côté droit, on a la représentation $\pi^{\GL}(\psi)$ qui est une induite ordinaire d'un caractère quadratique d'un sous-groupe de Levi. C'est donc bien une représentation faiblement unipotente. Il faut encore s'assurer que $\vert \nu\vert $ vu comme un caractère infinitésimal pour $G$ est inférieur à $\vert \nu_\psi\vert $ seulement si cela est vrai pour les transfert à $\GL(n^*)$. Le plus rapide est d'écrire la norme, en identifiant un caractère infinitésimal pour $G$ à $n$ nombres complexes (ici $n$ est le rang de $G$) et un caractère infinitésimal pour $\GL(n^*)$ à $n^*$ nombres complexes, le transfert envoie $\nu_n, \cdots, \nu_1$ sur $\nu_n, \cdots, \nu_1,\nu_0,-\nu_1, \cdots, -\nu_n$ où $\nu_0$ n'intervient que si $n^*=2n+1$ et vaut alors $0$. Et la norme (à un scalaire près) est la racine carré de $\sum_i \nu_i^2$ pour $G$ et de 2 fois cette somme pour $\GL(n^*)$.

\subsection{Description plus explicite des paquets d'Arthur\label{thprincipal}}
Le défaut du théorème de \ref{calcul} est que la translation est une opération dans le groupe de Grothendieck et n'est donc pas un foncteur exact. Le but de ce paragraphe est de montrer que dans \ref{calcul} on peut utiliser le foncteur de translation usuel, qui est donc un foncteur exact. Reprenons les notations de \ref{calcul} que l'on rappelle: on a un morphisme d'Arthur $\psi$ et un morphisme $\psi_+$ très régulier dominant $\psi$. D'où le caractère algébrique  $\lambda(\psi_+,\psi)$ . On note $E_+$ la représentation irréductible de $G$ de plus bas poids $-\lambda(\psi_+,\psi)$. 
On décompose $\psi$ comme en \ref{premiereforme}, $$\psi=\bigoplus_{i\in [1,v]}\rho_i \otimes R_{a_i}\bigoplus \psi_u.$$
On note encore $t_i$ les paramètres positifs des séries discrètes $\rho_i$. On considère les couples $(\mathfrak{q},L)$ de paires paraboliques avec $$L\simeq \times_{i\in [1,v]}U(p_i,q_i) \times G_0,$$
où pour tout $i\in [1,v]$, $p_i+q_i=a_i$ et $G_0$ est un groupe classique convenable. Et pour un tel $L$ on définit un caractère, $\lambda$, du produit des groupes unitaires en utilisant les paramètres $t_i$ pour $i\in [1,v]$ comme on l'a expliqué au début de cette section. On fixe aussi une représentation irréductible, $\sigma$, faiblement unipotente de $G_0$; elle n'est pas associée à $\psi_u$ mais à $\psi_u\otimes \chi$ où $\chi$ est un caractère quadratique qui dépend des $\rho_i$ mais cela n'a pas d'importance ici. 

En remplaçant $\psi$ par $\psi_+$ on remplace les $t_i$ par des demi-entiers $t'_i$, $\lambda$ par $\lambda'$ et si $\psi_+$ est très régulier en \cite{recent} 9.3, on montre que $\pi^A(\psi_+)$ est combinaison linéaire avec la propriété de multiplicité un des représentations $A_{\mathfrak{q}}(\lambda'\times \sigma)$ quand $\mathfrak{q}$ et $\sigma$ varient.  Cela ne sert que pour justifier l'hypothèse du théorème et surtout du corollaire ci-dessous.

\begin{theorem} Avec les notations précédentes, 
$$\left(A_{\mathfrak{q}}(\lambda_+ \times \sigma) \otimes E_+ \right)_{\nu_\psi}= A_{\mathfrak{q}}(\lambda \times \sigma).$$
\end {theorem}

\begin {cor}
On suppose  que $\pi^A(\psi_+)$ est une combinaison linéaire de repré\-sentations de $G$ de la forme $A_{\mathfrak{q}}(\lambda_+ \times \sigma)$ (cf. \cite{recent}) alors la même propriété est vraie pour $\pi^A(\psi)$ avec les mêmes coefficients (qui sont des caractères de $A(\psi)$) et $\pi^A(\psi)$ s'obtient en appliquant le foncteur de translation à $\pi^A(\psi_+)$ comme dans le théorème précédent.
\end {cor}
Disons tout de suite que le corollaire se déduit du théorème grâce à \ref{translation}.

Montrons le théorème.
On reprend le lemme 8.32 de \cite {KV}. On note $\mathfrak{u}$ le radical nilpotent de $\mathfrak{q}$. On note encore $L$ le sous-groupe de Levi de $\mathbb{q}$ défini sur $\mathbb{R}$ que l'on identifie à ses points réels et $L_1$ le groupe dérivé de $L$.
 Soit $\mu$ une combinaison linéaire de racines dans $\mathfrak{u}$. Et on suppose que $\mu$ est dominant pour $L$, c'est-à-dire que pour toute racine $\alpha$ dans $L$, on a $\langle \alpha,\mu\rangle\geq 0$. On note $E_\mu$ la représentation de dimension finie de $L$ de plus haut poids $\mu$, alors $$
\left( A_{\mathfrak{q}}((\lambda \times \sigma)\otimes E_\mu)\right)_\nu=0
$$
si $\mu\neq 0$.
En effet, on écrit le caractère infinitésimal $\nu$ comme le caractère infinitésimal associé au caractère du tore $\left(\lambda + \delta_{L,1}\right) + \nu_\sigma$, où $\delta_{L,1}$ est la demi-somme des racines de $L_1$. On décompose $\mu$ en $\mu_1 \times \mu_0$ suivant la décomposition de $L$ en deux facteurs, $\mu_1$ correspondant au produit des facteurs qui sont des groupes unitaires et $\mu_0$ correspond au facteur $G_0$. On note $E_{\mu_1}$ et $E_{\mu_0}$ les représentations de dimension finie de plus haut poids $\mu_1$ et $\mu_0$ respectivement.  La représentation $\lambda \otimes E_{\mu_1}$ a un caractère infinitésimal, c'est la représentation irréductible $E_{\mu_1}$ tordue par le caractère $\lambda$. C'est le caractère infinitésimal correspondant à $\lambda+\mu_1+\delta_{L_1}$. La représentation de $G_0$, $\sigma \otimes E_{\mu_0}$ est de longueur finie et ses sous-quotients irréductibles ont des caractères infinitésimaux de norme supérieure ou égale à  la norme de  $\nu_\sigma$. Montrons que la norme de $\lambda+\mu_1+\delta_{L_1}$ est strictement plus grande que la norme de $\lambda +\delta_{L_1}$ sauf si $\mu_1=0$. 
$$
 \mid\lambda+\mu_1+\delta_{L_1}\mid ^2= \mid \lambda +\delta_{L_1}\mid^2 +\mid \mu_1 \mid ^2 +2 \langle\lambda +\delta_{L_1}, \mu_1\rangle.
 $$
 Or $\langle\lambda,\mu_1\rangle
\geq 0 $ car $\mu_1$ est une combinaison linéaire de racines dans $\mathfrak{u}$. Et $\langle\delta_{L_1},\mu_1\rangle \geq 0$ car $\mu_1 $ est dominant pour $L_1$.

Ainsi la représentation $\left( A_{\mathfrak{q}}((\lambda \times \sigma)\otimes E_\mu)\right)$ est de longueur finie avec ses sous-quotients irréductibles ayant un caractère infinitésimal de norme strictement supérieure  à la norme de $\nu$ sauf éventuellement  si $\mu_1=0$. Si $\mu_1=0$, en fait on a $\mu=0$ car $\mu$ est une somme de racines dans $\mathfrak{u}$.

Venons en à la preuve du théorème en suivant maintenant \cite {KV} 8.35.
 On considère la représentation de dimension finie de $G$ de plus bas poids $-\lambda_+$. On la restreint à $(\mathfrak{q},L) $ et on la filtre  de sorte que les sous-quotients soient des représentations irréductibles de $L$ annulées par $\mathfrak{u}$. 
On remarque que les plus hauts poids de ces sous-quotients irréductibles sont de la forme $-\lambda_++\mu$ où $\mu$ vérifie exactement les conditions du paragraphe précédent. Ainsi la représentation
$$\left( A_{\mathfrak{q}}(\lambda \times \sigma)\right) \otimes E_{\lambda_+} \eqno(*)$$ est filtrée par des sous-représentations avec comme quotient des représentations de la forme
$$A_{\mathfrak{q}}\left((\lambda \times \sigma)\otimes E_\mu \right).\eqno(**)$$
 Donc la projection de la représentation (*) sur le caractère infinitésimal $\nu$ est aussi filtrée par des sous-représentations ayant comme quotient  les projections  des représentations (**) sur le caractère infinitésimal $\nu$. D'après ce qui précède toutes ces projections sont nulles sauf celles correspondant à $\mu=0$. Cela termine la preuve.

\end{document}